# Comment: Classifier Technology and the Illusion of Progress—Credit Scoring

**Ross W. Gayler**

These comments support Hand's argument for the lack of practical progress in classifier technology by pursuing them a little deeper in the specific context of credit scoring. Academic development of modeling techniques tends to ignore the role of the practitioner and the impact of business objectives. In credit scoring it can be seen that the nature of the task forces practitioners to adopt modeling strategies that positively favor simple techniques or, at least, limit the possible advantage of sophisticated techniques. The strategies adopted by credit scorers can be viewed as a heuristic approach to inference of the unobserved (and unobservable) distribution of possible data sets. The technical progress examined by Hand has been aimed toward better goodness of fit. However, technical progress toward a more principled basis for inferring the distribution of future problem data would be more likely to be adopted in practice.

## 1. CREDIT SCORING

I am approaching this commentary as a domain-specific consumer of statistical technology. My concern is credit scoring (the use of predictive statistical models to control operational decision-making in consumer finance). Classical credit scoring is applied at the point of application for a loan to predict the risk of default (nonpayment) and to make the decision whether to approve that application for credit. The total value of the loans made under the control of

*Ross W. Gayler is Honorary Associate, School of Communication, Arts and Critical Enquiry, La Trobe University, Melbourne, Australia and Senior Research and Development Consultant, Baycorp Advantage, Melbourne, Australia. Mailing address: 102 Through Road, Camberwell VIC 3124, Australia*

credit scoring is immense, and the value added to the economy by better decision-making because of credit scoring is correspondingly large. Thus, credit scoring is a domain where improved decision-making due to better predictive modeling would be valuable and technical progress would be expected.

Somewhat surprisingly, the statistical techniques currently used in credit scoring seem rather old-fashioned (often being simple regression models). This is not for lack of attempts to change the state of the art. New modeling techniques are regularly proposed for credit scoring (typically by academic researchers), but they are rarely adopted in practice. This lack of uptake cannot be blamed entirely on conservatism in the credit scoring community. The rewards of improvement are sufficiently high that once any lender adopts a better technique, there will be high competitive pressure for other lenders to do likewise. Rather, the continued use of simple predictive modeling techniques suggests that they have a practical advantage over more sophisticated techniques in credit scoring. Understanding the reasons for this advantage would be useful for the practice of applied predictive modeling in credit scoring and, more generally, might suggest productive avenues for the development of predictive modeling techniques to be applied in practical domains.

Professor Hand has worked extensively in credit scoring and it is likely that his experience in that domain motivated the writing of his paper, although his thesis, as stated, is not restricted to credit scoring. As a practitioner of credit scoring, I agree with the points he has raised. My aim here is to examine Hand's points a little further in the specific context of credit scoring, looking at the interaction of the technicalities of modeling with the demands imposed by the nature of the business task.

A brief description of the classical credit scoring problem is as follows. When credit is granted to consumers, some of the borrowers will default on their loans. The lender typically takes a loss on a defaulted loan. Ideally, a lender would predict which applicants would default and decline their applications for credit, thus avoiding the loss. The lender uses data available at the time of application to make that prediction and decision. The data may come from an application form, a credit bureau and the lender's own records if the applicant is an existing customer.

The potential predictors available at the time of application are not causally related to the outcome of default. Consequently, credit scoring models are correlational rather than causal. The outcome of default is not just dependent on the characteristics of the borrower, but also on external factors such as subsequent lender management actions and the state of the economy. Furthermore, the data are processed by the operational systems of lenders. These systems are constructed with the primary objective of carrying out the operational actions. Data collection and data quality issues that are relevant to statistical modeling are often an afterthought in system design (if they are considered at all). Consequently, the data quality is often not what would be desired, and data quality problems can be quite dynamic, because changes are made to the systems to accommodate short term operational needs. The data are noisy, and the quality of the noise is subject to drifts and jumps.

## 2. REGRESSION RATHER THAN CLASSIFICATION

Given that the occurrence of default is a binary outcome, it seems natural to treat credit scoring as a classification problem, and many academic papers have done so. Assuming a classification framework comes close to assuming that there is some ideal predictor space in which the outcome classes are perfectly separated. Even if such a predictor space does actually exist, it is not available to the credit scoring practitioner. The available predictors are not causally related to the outcome and some predictors (e.g., account management actions and changes in the economy) are not available at the time of the application because they occur subsequently. For problems such as this, as Hand notes more generally, "the Bayes error rate is high: meaning that no decision surface can separate the distributions of such problems very well" (Section 2.3). Given that the outcome classes cannot be separated, it may be better to adopt a regression framework for modeling and predict the probability of default conditional on the predictors.

However, in credit scoring there is an even more important consideration than the match between the theoretical form of the model and the true state of affairs. Lenders need to be able to control the rate at which loan applications are declined. This allows them to adjust workloads and to control the trade-off of profit against volume of business. A classification model yields predictions of "default" or "repay" which are mapped to decisions to "decline" or "accept" the loan application. Consequently, the decline rate is fixed by the predictions and the lender has no direct control of the decline rate from a classification model. This illustrates the point that credit scoring practitioners need to be mindful of the operational requirements of lending over and above goodness of fit and the theoretical form of models.

Hand's paper is written in terms of classifiers, but his arguments apply just as well to regression models used as classifiers. A regression model may be trivially converted to a classifier by having the predicted outcome be the probability of class membership and comparing it to a threshold. In fact, this is the standard form of credit scoring models. Conversely, some classification models can be converted to adequate regression models, but this is not generally true. A decision tree with two leaves will never make a good regression model. Consequently, even though classification models are not well suited to credit scoring, Hand's arguments do apply to credit scoring as it is practiced.

## 3. EQUIVALENCE OF MODELS AND DEGREES OF FREEDOM IN THE MODELER

Hand observed that "a tremendous variety of algorithms and models has been developed for the construction of such [classification] rules" (Section 1). Different algorithms have different representational biases and a different bias/variance trade-off. For a fixed set of predictors we would expect different algorithms to generate different approximations to the outcome. However, in credit scoring the set of predictors is not fixed. The model developer is free to generate new derived variables in the data set and will generally do so to accommodate the particular representational bias of the modeling technique used. For example, decision tree induction and projection pursuit regression are able to automatically model interactions in the data, whereas regression works only with the predictors it is given and does not create interactive combinations. The credit scoring modeler using regression would construct interaction predictors if they were thought necessary.

The objective of every modeling technique is to approximate the data. Thus, in the limit (and the hands of a skilled modeler), every modeling technique should end up in agreement because they are all approximating the same data. However, the effort required to achieve that degree of approximation may vary greatly between techniques. Even for techniques that require the same effort to achieve a given accuracy of approximation, the models may differ in other properties that are operationally important to the lender.

It is also worth recalling Hand's comment about the high Bayes error rate (Section 2.3). When the ratio of variance accounted for by the response surface is low compared to the error about the response surface (as it is in credit scoring), it becomes harder to distinguish between different representational biases. Thus we would not expect the differences between different modeling techniques to be readily observable.

The impact of the skilled modeler warrants some further investigation. Effectively, the modeler supplies extra degrees of freedom in addition to those supplied by the modeling technique. The natural consequence of this is to reduce the difference between techniques in terms of goodness of fit. Rather than compare modeling techniques in terms of predictive power, it would be more useful to look at the effort required of the modeler to achieve a given goodness of fit and other properties of the models that are of operational relevance to the lender.

## 4. MAIN EFFECTS AND INTERACTIONS

In Section 2.3, Hand mentions "examples of artificial data which simple models cannot separate (e.g., intertwined spirals or checkerboard patterns)," noting that "such data sets are exceedingly rare in real life [and] it is common to find that the centroids of the predictor variable distributions of the classes are different." This is a claim that problems which can be modeled only as interactions of the variables (with no observable main effects) are rare. This may well be true in general because of the improbability of interactions exactly canceling out to leave no main effects. However, in credit scoring it is also true for domain-specific reasons. The inclusion of each predictor in a decision-making system has to be justified (operationally and legally). It is much easier to argue for the inclusion of a predictor if the argument can be made for that predictor in isolation. Conversely, it is harder to argue for the inclusion of a predictor if it can be shown to add value only in the context of other predictors.

Furthermore, credit scoring practitioners are very concerned with the stability over time of their models. Some credit scoring models are used for years before being replaced. Therefore, it is important to ensure that the predictive relationships on which the model is based are stable over time. Credit scoring practitioners tend to believe that main effects are more stable than interactions (all other things being equal). When interactions are included as predictors, it is generally because the modeler has a prior belief that the interaction reflects some stable mechanism in the world. An otherwise unmotivated interaction that is discovered by an automated search procedure is unlikely to be included in a predictive model or, if it is included, to have its influence intentionally limited relative to the main effects. The effect of these selection biases is to ensure that credit scorers prefer simpler models based on main effects.

## 5. SENSITIVITY TO ARBITRARY MODELING DECISIONS

Hand notes that when constructing classification rules, "various ... assumptions and choices are often made which may not be appropriate" (Section 1) and even when they are entirely appropriate, the choices may be somewhat arbitrary. He gives the example of typically defining "a customer as 'defaulting' if they fall three months in arrears with repayments ... [while] [i]t is entirely reasonable that alternative definitions (e.g., four months in arrears) might be more useful if economic conditions were to change" (Section 4.2). Credit scoring necessarily involves many detailed decisions concerning the modeling process. Many of these decisions involve compromises and trade-offs, with no obviously correct answer. While the experienced credit scorer would have arguments for the specific decisions made, it would be a bold modeler who would argue that the decisions taken were uniquely and obviously correct. Thus, there is an element of arbitrariness in the modeling process.

It is possible to conceive of a space of feasible modeling decisions. Similar sets of decisions are nearby in that space. A small change in the modeling decisions would generally lead to a small change in the models. However, the possibility exists that a small change in modeling decisions may lead to a large change in the models that arise from them. This would be very unsatisfactory in credit scoring because the results of the modeling would be strongly dependent on arbitrary modeling choices. Therefore, credit scorers tend to restrict their attention to regions of the modeling decision space where the gradient of models with respect to modeling decisions is low. In these regions, all the models generated as a result of the different modeling choices would yield similar results. If a new modeling technique yielded markedly different results, it would be unlikely to be favored by credit scorers unless it was surrounded by a region of other models yielding similar results. It would be more difficult for the modeler to argue for the correctness of the unique results given that the choice of modeling technique might be regarded as arbitrary.

## 6. DEVELOPMENT DATA NOT REPRESENTATIVE OF OPERATION

Hand points out "that in many ... real classification problems the data points in the design set are not ... randomly drawn from the same distribution as the data points to which the classifier will be applied" (Section 1). Furthermore, any design set represents "merely a single ... problem drawn from a notional distribution of problems" (Section 1). Later he notes that "a fundamental assumption of the classical paradigm is that the various distributions involved do not change over time ... [although this assumption] is unrealistic in most commercial applications, concerned with human behaviour" (Section 3.1). His concern here is with population drift. This would not be a problem if the predictive model were the "true" model, but as Hand states "it would be a brave person who could confidently assert that [this] held" (Section 3.2).

Population drift is a particular concern in credit scoring. Loans which default do so over an extended period after the loan has been granted. Consequently, an extended outcome period (typically at least one year) is required to allow a reasonable proportion of loans to default. To this must be added time to accumulate enough applications to provide a reasonable number of observations for modeling and to allow for seasonal variation in the applicant population. Allowing time for data preparation, data modeling and implementation of the models into the operational system, it is common for the oldest data on which a model is based to be three years old when the model is first switched on. Then the model may be in use for some while (three years is common, and more than five years not unknown). Even if the applicant population distribution is stationary, the data collecting process is subject to random jumps, because lenders may change their systems and procedures at any time. Thus, a large part of the value added by credit scoring practitioners comes from anticipating possible future shifts in the data distribution and designing the models to be relatively insensitive to such shifts. This can be seen as another aspect of attempting to reduce the sensitivity of the models to arbitrary features of the specific design set (in this case, characteristics of the data that just happen to hold at the time the data are collected).

The expertise of the credit scoring modeler can be thought of as applying a bias to the modeling techniques to move the models toward the notional distribution of problems. For example, Hand discusses the application of a tree model and linear discriminant analysis (as competing techniques) to consumer credit data, and points out that because the design set is always retrospective, the population may have drifted by the time the model is built and "reduced any advantage that the more sophisticated tree model may have" (Section 3.1). A tree model fits better than linear discriminant analysis, but degrades more rapidly. There is the possibility that the tree model may actually become worse than the linear discriminant model with the passage of time. Rather than view the techniques as competing, a credit scorer might model the data with linear discriminant analysis and then build a tree model of the residuals. This hybrid model puts a bound on deterioration by predicting the majority of the outcome variance using the more stable modeling technique.

## 7. FREEDOM VIA THE FLAT MAXIMUM EFFECT

Hand mentions the flat maximum effect in the context of explaining that a reasonable fraction of the maximum attainable predictive power can be obtained from an equally weighted combination of predictors (Section 2.4). The existence of the flat maximum effect is a great advantage in credit scoring. It implies that there may be many alternative models with similar goodness of fit. This provides the credit scoring modeler the opportunity to choose between those models on some basis other than goodness of fit (e.g., susceptibility to population drift or ability to finely control the decline rate). The freedom this confers is so valuable that credit scoring modelers prefer to choose predictors that make the flat maximum effect more likely to exist. This is the case where there is a conditional monotone relationship between each of the predictors and the outcome (which also happens to be the circumstances under which a simple linear combination is likely to perform well).

## 8. VALUE ADD AND MODELING TECHNIQUES

In credit scoring, much of the value added by modelers is not via goodness of fit to the development sample, but by anticipation of possible changes in the operational systems and data. This can be viewed as a problem of trying to infer the unobserved distribution of possible development data sets. Credit scorers attempt to achieve this by biasing their models toward simple models and techniques. These models are not only more likely to generalize across potential data sets, but also, as Hand points out, to yield most of the predictive power of more complex models. More complex models of the current data set are unlikely to be attractive to credit scorers. However, techniques that provide a more principled basis for generalizing to the distribution of possible data sets would be welcome.